\documentclass[11pt,a4paper]{article}
\usepackage[T1]{fontenc}
\usepackage{amsmath,amssymb,amsthm}
\usepackage{tikz}
\usetikzlibrary{arrows.meta}
\usepackage{cite}
\usepackage[hidelinks,bookmarksnumbered,bookmarksopen,hypertexnames=false]{hyperref}
\usepackage{doi}

\date{}

\makeatletter
\def\namedlabel#1#2{\begingroup
  \def\@currentlabel{#2}\def\@currentHref{#1}\Hy@raisedlink{\hypertarget{#1}{}}#2\label{#1}\endgroup
}
\makeatother

\newtheorem{theorem}{Theorem}[section]
\newtheorem{lemma}[theorem]{Lemma}
\theoremstyle{definition}
\newtheorem{definition}[theorem]{Definition}
\newtheorem{example}[theorem]{Example}
\newtheorem{remark}[theorem]{Remark}
\numberwithin{equation}{section}

\newcommand{\R}{\mathbb{R}}
\DeclareMathOperator*{\essinf}{ess\,inf}
\DeclareMathOperator*{\esssup}{ess\,sup}
\let\div\undefined\DeclareMathOperator{\div}{div}

\title{
Multiplicity of positive solutions for degenerate nonlocal $p$-Laplacian problems without monotonicity of the normalized reaction
}
\author{Giuseppe Failla$^{1}$, Leszek Gasi\'nski$^{2}$, Jo\~ao R. Santos J\'unior$^{3}$}

\begin{document}
\maketitle
\noindent$^{1}$ Department of Mathematics and Computer Sciences, University of Palermo, Via Archirafi 34, 90123 Palermo, Italy.\\
E-mail: giuseppe.failla02@unipa.it\\
$^{2}$ Department of Mathematics, University of the National Education Commission, ul. Podchor\k a\.zych 2, 30-084 Krak\'ow, Poland.\\
E-mail: leszek.gasinski@uken.krakow.pl\\
$^{3}$ Instituto de Ci\^{e}ncias Exatas e Naturais, Universidade Federal do Par\'a, Avenida Augusto Corr\^{e}a 01, 66075-110 Bel\'em, PA, Brazil.\\
E-mail: joaojunior@ufpa.br

\begin{abstract}
We study degenerate nonlocal Dirichlet problems driven by the
$p$-Laplacian, with diffusion depending on weighted sums of integral powers
of the solution. Under explicit sign, spectral and size conditions,
we obtain two positive solutions for each bounded positivity interval of
the diffusion coefficient with zero endpoint values. Neither monotonicity
of the normalized reaction $f(x,t)/t^{p-1}$ nor uniqueness of the frozen
problem is required. The proof combines Schauder's theorem on a
product space of a scalar parameter and a function
with subsolution bounds and an energy estimate.
The multiplicity result also persists under uniformly small
continuous perturbations of the nonlocal functional, including
perturbations that do not preserve order.
\end{abstract}



\section{Introduction}
We consider the nonlocal Dirichlet problem
\begin{equation}\tag{P}\label{Problem1}
\left\{
\begin{array}{ll}
 -m\bigl(\mathcal H(u)\bigr)\Delta_pu=f(x,u)&\hbox{in }\Omega,\\
 u>0&\hbox{in }\Omega,\\
 u=0&\hbox{on }\partial\Omega,
\end{array}\right.
\end{equation}
where $\Omega\subset\R^N$ is a bounded domain with smooth boundary,
$1<p<+\infty$, and $\Delta_pu=\div(|\nabla u|^{p-2}\nabla u)$.
The nonlocal functional is
\begin{equation}\label{eq:H-functional}
 \mathcal H(u)=\sum_{i=1}^M\int_\Omega\mu_i(x)|u|^{q_i}\,dx,
 \qquad u\in L^\infty(\Omega),
\end{equation}
with $M\geq1$, $q_i\geq1$ and nonnegative bounded measurable weights
$\mu_i$. We seek two positive solutions for each bounded interval on
which $m$ is positive and vanishes at both endpoints.

Equations of this type describe stationary states of models in which a
spatially uniform diffusion rate depends on an integral of the state.
For example, \eqref{Problem1} is the stationary equation associated with
\[
 \partial_tu-\div\bigl(m(\mathcal H(u))|\nabla u|^{p-2}\nabla u\bigr)
 =f(x,u).
\]
This interpretation justifies the terms \emph{diffusion
coefficient} for $m$ and \emph{reaction} for $f$ used below. The nonlocal
dependence is in the scalar coefficient, while
the differential operator is local. Classical motivations include
Carrier's string models \cite{Car1,Car2} and population dispersal models
\cite{ChiRod,ChiLov2}. More recent work studies how spatial heterogeneity
and degeneracy of the coefficient affect positive stationary states
\cite{FigFigMorSua,DelMorSanSua}, as well as nonlocal evolution equations
\cite{CarChi}. In the present setting, weights such as
$\mu_i=\chi_{\Omega_i}$ allow the coefficient to depend on measurements
over prescribed regions. These regions may overlap and need not cover
$\Omega$. Our result concerns the multiplicity of stationary states and
does not address their dynamical stability.

Degeneracy of the nonlocal coefficient may change the structure
of the positive solution set. For a logistic reaction, Delgado et al.\
\cite{DelMorSanSua} show that the continuum of positive solutions of
the local problem is replaced by two disjoint branches in the nonlocal
problem. 
For degenerate coefficients, multiplicity
results with two positive
solutions per positivity interval were obtained for the Laplacian in
\cite{GasSan1,GasSan2} and for the $p$-Laplacian in
\cite{CanGasLivSan}. Related extensions include variable exponent,
Neumann and double phase problems \cite{CanFaiLiv,CreFaiVas,FaiGas}.
In \cite{CanGasLivSan}, the nonlocal term is $\int_\Omega|u|^q\,dx$
and the autonomous reaction satisfies a strict decrease condition on
$f(t)/t^{p-1}$. This condition yields uniqueness for the nonlinear
problem obtained by fixing the nonlocal argument. Its solution
$u_\alpha$ then defines a continuous scalar map
$\alpha\mapsto\int_\Omega u_\alpha^q\,dx$, whose fixed points give the
required solutions.

Gasi\'nski, Santos J\'unior and Siciliano \cite{GasSanSic} treat the
semilinear equation
$-\Delta u=F(u,\int_\Omega|u|^q\,dx)$, allowing singular dependence of
$F$ on its second argument. Assumption (f2) of that paper
requires $s\mapsto F(s,\alpha)/s$ to be decreasing. Its Proposition 2.1
gives a unique positive solution for each admissible $\alpha$, and its
Propositions 3.1 and 3.2 use the resulting continuous scalar map to obtain
$2K$ solutions.
Here the reaction may depend measurably on $x$, the operator is the
$p$-Laplacian, and no decrease condition is imposed on the normalized
reaction. We retain the same multiplicity pattern through a different
fixed point construction. The two results have different scopes, since
we do not treat the general singular dependence of $F$ on its second
argument allowed in \cite{GasSanSic}.

Our main contribution is to obtain this multiplicity conclusion without
assuming monotonicity of
\begin{equation}\label{eq:intro-monotonicity}
 (0,t_*)\ni t\longmapsto\frac{f(x,t)}{t^{p-1}}.
\end{equation}
We refer to this quotient as the \emph{normalized reaction}.
Fixing the nonlocal argument at a value $\alpha$ produces the auxiliary
Dirichlet problem
\begin{equation}\label{eq:intro-frozen}
 -m(\alpha)\Delta_pu=f(x,u)\quad\hbox{in }\Omega,
 \qquad u=0\quad\hbox{on }\partial\Omega,
\end{equation}
which we call the \emph{frozen problem}. If the normalized
reaction is not decreasing, uniqueness of positive solutions of
\eqref{eq:intro-frozen} does not follow from the principles of Br\'ezis
and Oswald \cite{BreOsw} or D\'iaz and Sa\'a \cite{DiaSaa}, and no
continuous selection $\alpha\mapsto u_\alpha$ of solutions is available
from those arguments. We instead use an operator $\mathcal T(\alpha,v)$
defined by solving the Dirichlet problem obtained from
\eqref{eq:intro-frozen} when the right-hand side is evaluated at a given
function $v$ and suitably truncated. Schauder's theorem is applied to a
map defined on the product of a parameter interval and a ball in
$C^1_0(\overline\Omega)$. A fixed point $(\alpha,u)$ of this map
satisfies $u=\mathcal T(\alpha,u)$, so that $u$ solves the frozen problem
at $\alpha$. Bounds at the two ends of the parameter interval then
enforce the additional identity $\alpha=\mathcal H(u)$, which turns $u$
into a solution of \eqref{Problem1}.
Lemma~\ref{lem:abstract-constraint-fp} records this elementary consequence
of Schauder's theorem.

Subsolution and supersolution methods are already used for nonlocal
elliptic equations, for example in \cite{AlvCov}, where they are combined
with pseudomonotone operator theory. In our proof, the lower
barrier is a subsolution of the frozen problem which varies continuously
with $\alpha$. Near the endpoints of each positivity interval it is chosen
sufficiently large to make $\mathcal H(u)>\alpha$ for every fixed point
$u$ of the truncated operator at such a value of $\alpha$, whereas an
energy estimate gives $\mathcal H(u)<\alpha$ at a suitable interior point
$\alpha=\tau_k$. Both estimates apply to all fixed points at the specified
endpoints, and not only to a selected one. This is what makes the
argument insensitive to the possible nonuniqueness for the frozen
problem. Applying the construction on each
side of the interior point gives two distinct solutions. The same
argument accommodates the weighted sum in \eqref{eq:H-functional}, with
different exponents and weights that need not be bounded away from zero.

Suitably constructed subsolutions allow us to formulate the size
condition in terms of the total weighted mass
$M_*=\sum_i t_*^{q_i}\int_\Omega\mu_i$. A continuous family connects
these barriers to a small multiple of the first eigenfunction, for every
$p>1$. Theorem~\ref{thm:perturbed-functional} further shows that the
multiplicity conclusion persists under continuous perturbations of
$\mathcal H$ that are uniformly small on $0\leq u\leq t_*$. The perturbed
functional need not preserve order. Order preservation of the reference
functional supplies the strict endpoint estimates used in the proof.
Example~\ref{ex:nonorder-functional} gives a perturbation depending on an
additional spatial integral and proves that the resulting coefficient
cannot be expressed as a function of the reference value $\mathcal H(u)$
alone.

\begin{definition}\label{def:weak-solution}
  A function $u\in W^{1,p}_0(\Omega)\cap L^\infty(\Omega)$ is a
  \emph{(weak) positive solution} of \eqref{Problem1} if $u>0$ a.e. in
  $\Omega$ and
\[
  m\big(\mathcal H(u)\big)
  \int_\Omega|\nabla u|^{p-2}\nabla u\cdot\nabla\varphi\,dx
  =\int_\Omega f(x,u)\varphi\,dx
  \qquad\hbox{for all }\varphi\in W^{1,p}_0(\Omega).
\]
\end{definition}

We impose the following assumptions on the weights.
\begin{itemize}
  \item[\namedlabel{itm:Hmu}{$(H_\mu)$}] For every $i\in\{1,\ldots,M\}$,
  $\mu_i\in L^\infty(\Omega)$, $\mu_i\geq0$ a.e. in $\Omega$, and
  there exists at least one $i\in\{1,\ldots,M\}$ such that
  $\mu_i\not\equiv0$.
\end{itemize}

The coefficient satisfies the following condition.
\begin{itemize}
  \item[\namedlabel{itm:Hm}{$(H_m)$}] $m:[0,+\infty)\to\R$
  is continuous and there exist numbers
  \[
    0\leq\alpha_1<\alpha_2\leq\alpha_3
    <\alpha_4\leq\cdots\leq\alpha_{2K-1}<\alpha_{2K},\qquad K\geq1,
\]
  such that
  \[
    m>0\quad\hbox{in }(\alpha_{2k-1},\alpha_{2k}),\qquad
    m(\alpha_{2k-1})=m(\alpha_{2k})=0
\]
  for all $k\in\{1,\ldots,K\}$.
\end{itemize}
The positivity intervals in \ref{itm:Hm} are pairwise disjoint. No sign
condition is imposed on $m$ outside them. Every solution furnished by Theorem~\ref{thm:main}
satisfies $m(\mathcal H(u))>0$.

We assume the following properties of the reaction.
\begin{itemize}
 \item[\namedlabel{itm:Hf}{$(H_f)$}]
 The function $f:\Omega\times\R\to\R$ is Carath\'eodory and satisfies
 the following conditions for some $t_*>0$.
 \begin{itemize}
  \item[(i)]
  $f(x,t_*)=0$ for a.e. $x\in\Omega$.
  \item[(ii)] For every compact interval $J\subset(0,t_*)$,
  \[
   \inf_{t\in J}\ \essinf_{x\in\Omega}f(x,t)>0.
  \]
  \item[(iii)] There exists $T>0$ such that
  \[
   \esssup_{(x,t)\in\Omega\times[0,t_*]} f(x,t)\leq T.
  \]
 \end{itemize}
\end{itemize}

Let $\lambda_1$ be the first Dirichlet eigenvalue of $-\Delta_p$, and let
$e_1$ be its positive eigenfunction normalized by $\Vert e_1\Vert_\infty=1$.
We require three compatibility conditions.
\begin{itemize}
  \item[\namedlabel{itm:Hlambda}{$(H_\lambda)$}]
  $\displaystyle
  \liminf_{t\to0^+}
  \essinf_{x\in\Omega}
  \frac{f(x,t)}{t^{p-1}}
  >
  \lambda_1
  \max_{k\in\{1,\ldots,K\}}
  \max_{\alpha
  \in[\alpha_{2k-1},\alpha_{2k}]}m(\alpha)$.
  \item[\namedlabel{itm:He}{$(H_M)$}]
  $\displaystyle
  \alpha_{2K}<M_*:=\sum_{i=1}^{M}t_*^{q_i}
  \int_\Omega \mu_i(x)\,dx=\mathcal H(t_*)$.
  \item[\namedlabel{itm:Htau}{$(H_\tau)$}] 
  For each $k\in\{1,\ldots,K\}$ there exists 
  $\tau_k\in(\alpha_{2k-1},\alpha_{2k})$ such that
\[
    \tau_k^p m(\tau_k)>
    \frac{|\Omega|}{\lambda_1}
    \left(\sum_{i=1}^{M}t_*^{q_i-1}
    \Vert\mu_i\Vert_{p'}\right)^p E,
\]
  where $E:=\esssup\limits_{(x,t)\in\Omega\times[0,t_*]}[t f(x,t)]$.
\end{itemize}
Here $p'=p/(p-1)$, $\Vert\cdot\Vert_{p'}$ denotes the norm in
$L^{p'}(\Omega)$, and $|\Omega|$ is the Lebesgue measure of $\Omega$.
Assumptions \ref{itm:Hf} and \ref{itm:Hmu} imply $0<E<+\infty$ and
$\Vert\mu_i\Vert_{p'}<+\infty$.

We now state the first main result.

\begin{theorem}\label{thm:main}
Assume \ref{itm:Hmu}, \ref{itm:Hm}, \ref{itm:Hf},
\ref{itm:Hlambda}, \ref{itm:He} and \ref{itm:Htau}.
Then problem \eqref{Problem1} has at least $2K$ pairwise distinct positive
solutions in $C^1_0(\overline\Omega)$. More precisely, for each
$k\in\{1,\ldots,K\}$ there are solutions $u_{k,1}$ and $u_{k,2}$ such that
$0<u_{k,j}\leq t_*$ in $\Omega$ for $j=1,2$, and
\begin{equation}\label{eq:main-order}
 \alpha_{2k-1}<\mathcal H(u_{k,1})<\tau_k
 <\mathcal H(u_{k,2})<\alpha_{2k}.
\end{equation}
\end{theorem}

The ordering in \eqref{eq:main-order} concerns the values of
$\mathcal H$, and does not assert a pointwise ordering of the solutions.
The role of the
assumptions is as follows. Assumption \ref{itm:Hlambda} provides a small
positive subsolution $\theta_0e_1$ of the frozen problem
\eqref{eq:intro-frozen} for every $\alpha$ in the positivity intervals of
$m$. Assumptions \ref{itm:Hlambda}, \ref{itm:Hf}(ii) and \ref{itm:He} provide larger
subsolutions for $\alpha$ near the zeros of $m$.
These give the inequality
$\mathcal H(u)>\alpha$ for the fixed points at the outer endpoints of the
parameter intervals used in the proof. Assumption \ref{itm:Htau} yields
the reverse inequality $\mathcal H(u)<\tau_k$ for every solution $u$ of
the frozen problem at $\alpha=\tau_k$ with $0\leq u\leq t_*$. No
condition on $f$ outside $[0,t_*]$ is used, and
the exponents $q_i$ are not restricted by a Sobolev critical exponent.
Example~\ref{ex:nonmonotone-nonlinearity} verifies these assumptions for
a normalized reaction that is neither nondecreasing nor nonincreasing.

\begin{remark}\label{rem:chi-special-case}
For measurable sets $\Omega_i\subset\Omega$ of positive measure, the
choice $\mu_i=\chi_{\Omega_i}$ gives
\[
 \mathcal H(u)=\sum_{i=1}^M\int_{\Omega_i}|u|^{q_i}\,dx
 =\sum_{i=1}^M\Vert u\Vert_{q_i,\Omega_i}^{q_i}.
\]
No regularity, disjointness or covering condition is required of these sets.
\end{remark}

The multiplicity conclusion persists for a class of perturbations that
need not preserve order. Let $X=C^1_0(\overline\Omega)$ and set
\[
 \mathcal D=\{u\in X:0\leq u\leq t_*\text{ in }\Omega\}.
\]

\begin{theorem}\label{thm:perturbed-functional}
Assume the hypotheses of Theorem~\ref{thm:main}.
There exists $\varepsilon_0>0$ such that, if 
$\widetilde{\mathcal H}:X\to[0,\infty)$ is continuous and satisfies
\begin{equation}\label{eq:perturb-smallness}
 \sup_{u\in\mathcal D}
 |\widetilde{\mathcal H}(u)-\mathcal H(u)|\leq\varepsilon_0,
\end{equation}
then the problem
\begin{equation}\label{eq:perturbed-problem}
 \left\{\begin{array}{ll}
 -m(\widetilde{\mathcal H}(u))\Delta_pu=f(x,u)&\text{in }\Omega,\\
 u>0&\text{in }\Omega,\\
 u=0&\text{on }\partial\Omega
 \end{array}\right.
\end{equation}
has at least $2K$ pairwise distinct positive solutions
$u_{k,j}\in X$, $k\in\{1,\ldots,K\}$, $j\in\{1,2\}$, with
$0<u_{k,j}\leq t_*$ in $\Omega$, and
\begin{equation}\label{eq:perturbed-order}
 \alpha_{2k-1}<\widetilde{\mathcal H}(u_{k,1})<\tau_k
 <\widetilde{\mathcal H}(u_{k,2})<\alpha_{2k}
 \quad(k=1,\ldots,K).
\end{equation}
\end{theorem}

\begin{remark}\label{rem:perturb-scope}
No order-preserving property is required of $\widetilde{\mathcal H}$.
The smallness assumption controls only the values of the functional on
$\mathcal D$, not its derivatives or any Lipschitz constant. Although
$\mathcal D$ is not bounded in the $C^1$ norm, the reference values satisfy
$0\leq\mathcal H(u)\leq M_*$ there. No perturbation bound outside
$\mathcal D$ is needed. Order preservation of the reference functional is
used to establish the strict endpoint margins, which then persist for
the actual functional. In particular, this theorem does not replace the
monotonicity of the differential operator $-\Delta_p$, and it does not
claim the same conclusion for arbitrary continuous nonlocal functionals.
The conclusion concerns persistence of multiplicity and
does not assert continuous dependence of individual solutions on the
perturbation. The proof shows that $\varepsilon_0$ may be any number in
$(0,d)$, where $d>0$ is the explicit margin defined in
\eqref{eq:perturb-margin}.
\end{remark}

Section~\ref{SecPrel} contains the fixed point lemma and the
subsolution construction. Section~\ref{SecAux} develops the auxiliary
problems and endpoint estimates. Section~\ref{SecMain} contains the
proofs of Theorems~\ref{thm:main} and~\ref{thm:perturbed-functional}.
Section~\ref{SecExample} gives three examples illustrating the main results.

\section{Preliminaries and auxiliary tools}\label{SecPrel}
For a measurable set $D\subseteq\Omega$ of positive measure, we use the
standard Lebesgue spaces $L^r(D)$ and norms $\Vert\cdot\Vert_{r,D}$.
When $D=\Omega$, the subscript $D$ is omitted. The space
$W^{1,p}_0(\Omega)$ is endowed with the norm
$\Vert u\Vert=\Vert\nabla u\Vert_p$, and its dual is denoted by
$W^{-1,p'}(\Omega)$. We write $s^+=\max\{s,0\}$ and
$s^-=\max\{-s,0\}$, and denote the characteristic function of a
measurable set $D$ by $\chi_D$. We also use
\[
 C^1_0(\overline\Omega)=\{u\in C^1(\overline\Omega):u=0\text{ on }\partial\Omega\},
 \qquad
 \Vert u\Vert_{C^1(\overline\Omega)}=\Vert u\Vert_\infty+\Vert\nabla u\Vert_\infty.
\]
For $\nu\in(0,1)$, the embedding
$C^{1,\nu}(\overline\Omega)\hookrightarrow C^1(\overline\Omega)$ is compact.

The operator $-\Delta_p:W^{1,p}_0(\Omega)\to W^{-1,p'}(\Omega)$ is
defined by
\[
 \langle-\Delta_pu,v\rangle
 =\int_\Omega|\nabla u|^{p-2}\nabla u\cdot\nabla v\,dx.
\]
It is continuous, bounded on bounded sets, coercive and strictly monotone
\cite{MotMotPap}. This operator monotonicity is distinct from the
decrease condition on the normalized reaction discussed in the
Introduction. The first eigenvalue satisfies
\begin{equation}\label{eq:lambda1}
 \lambda_1=\inf_{u\in W^{1,p}_0(\Omega)\setminus\{0\}}
 \frac{\int_\Omega|\nabla u|^p\,dx}{\int_\Omega|u|^p\,dx}>0.
\end{equation}
The corresponding normalized positive eigenfunction satisfies
\[
 -\Delta_pe_1=\lambda_1e_1^{p-1},\qquad
 e_1\in C^1_0(\overline\Omega),\qquad 0<e_1\leq1\text{ in }\Omega.
\]
These properties follow from the spectral theory of the
Dirichlet $p$-Laplacian, together with the regularity theory and the
strong maximum principle, see \cite[Chapters 8 and 9]{MotMotPap}.

We shall repeatedly use the following consequence of the
Carath\'eodory property of $f$. Suppose that a nonstrict inequality between two
functions of $(x,t)$ which are continuous in $t$ holds, for each fixed
$t$ in an interval $J$, for a.e. $x\in\Omega$. Then it holds for every
$t\in J$ outside a single null set of $x$. 
It suffices to take the union
of the exceptional null sets corresponding to a countable dense subset of
$J$ and to use continuity in $t$. 
If $J$ has nonempty interior,
a nonstrict inequality valid almost everywhere on $\Omega\times J$ has
the same property, by Fubini's theorem and continuity in $t$. 
In this way, \ref{itm:Hf}(iii) and the definition
of $E$ in \ref{itm:Htau} give
\begin{equation}\label{eq:f-upper-bound-uniform}
 0\leq f(x,t)\leq T
 \quad\hbox{for a.e. }x\in\Omega\hbox{ and every }t\in[0,t_*],
\end{equation}
and
\begin{equation}\label{eq:tf-upper-bound-uniform}
 t f(x,t)\leq E
 \quad\hbox{for a.e. }x\in\Omega\hbox{ and every }t\in[0,t_*],
\end{equation}
where the lower bound in \eqref{eq:f-upper-bound-uniform}
follows from \ref{itm:Hf}(ii), applied on a countable exhaustion of
$(0,t_*)$ by compact intervals, together with \ref{itm:Hf}(i) and
the continuity of $f(x,\cdot)$ at $t=0$. The same argument gives a
uniform positive lower bound for $f$ on $\Omega\times J$ for each compact
interval $J\subset(0,t_*)$. All these estimates remain valid when $t$ is
replaced by a measurable function with values in the relevant interval.

We use the following form of Schauder's theorem
\cite[Theorem 4.7]{MotMotPap}.

\begin{theorem}[Schauder fixed point theorem]\label{thm:schauder}
  Let $X$ be a Banach space, let $C\subset X$ be nonempty, closed, 
  bounded, and convex, and let $\mathcal K:C\to C$ 
  be continuous and compact, 
  in the sense that $\mathcal K(C)$ is relatively compact in $X$. 
  Then $\mathcal K$ has a fixed point.
\end{theorem}

\begin{lemma}[Fixed point principle for a scalar constraint]
\label{lem:abstract-constraint-fp}
Let $X$ be a Banach space, let $B\subset X$ be nonempty, closed,
bounded and convex, and let $\Phi:X\to\R$ be continuous.
Let $a<b$ and let $\mathcal T:[a,b]\times B\to B$ be continuous with
relatively compact image. Suppose that one of the following
two alternatives holds, the implications being required for every
$u\in B$.
\begin{enumerate}
 \item[(a)]
 \[
  u=\mathcal T(a,u)\ \Longrightarrow\ \Phi(u)>a,
  \qquad
  u=\mathcal T(b,u)\ \Longrightarrow\ \Phi(u)<b.
 \]
 \item[(b)]
 \[
  u=\mathcal T(a,u)\ \Longrightarrow\ \Phi(u)<a,
  \qquad
  u=\mathcal T(b,u)\ \Longrightarrow\ \Phi(u)>b.
 \]
\end{enumerate}
Then there exists $(\alpha,u)\in(a,b)\times B$ such that
\[
 u=\mathcal T(\alpha,u),\qquad \alpha=\Phi(u).
\]
\end{lemma}

\begin{proof}
Let $\Pi(s)=\min\{b,\max\{a,s\}\}$, and set $\sigma=1$ in case (a)
and $\sigma=-1$ in case (b). Define
\[
 \mathcal K(\alpha,v)=
 \left(\Pi\bigl(\alpha+\sigma[\Phi(\mathcal T(\alpha,v))-\alpha]\bigr),
       \mathcal T(\alpha,v)\right).
\]
This is a continuous map from $[a,b]\times B$, a nonempty,
closed, bounded and convex subset of the Banach space $\R\times X$, into
itself. Its image is
contained in the compact set
$[a,b]\times\overline{\mathcal T([a,b]\times B)}$, so $\mathcal K$ is a compact map.
Theorem~\ref{thm:schauder} gives a fixed point $(\alpha,u)$,
so that $u=\mathcal T(\alpha,u)$ and
$\alpha=\Pi\bigl(\alpha+\sigma[\Phi(u)-\alpha]\bigr)$. At $\alpha=a$ the
endpoint assumption gives $\sigma[\Phi(u)-a]>0$, so the first component
of $\mathcal K(a,u)$ is greater than $a$. At $\alpha=b$ it gives
$\sigma[\Phi(u)-b]<0$, so the first component is less than $b$.
Both possibilities contradict the fixed point identity. Hence
$\alpha\in(a,b)$. Since $\Pi(s)=\alpha$ with $\alpha\in(a,b)$
forces $s=\alpha$, the projection identity reduces to
\[
 \alpha=\alpha+\sigma[\Phi(u)-\alpha].
\]
Since $\sigma\neq0$, we obtain $\alpha=\Phi(u)$.
\end{proof}

The next two lemmas provide subsolutions for the frozen problem.

\begin{lemma}\label{lem:small-barrier}
  Assume that \ref{itm:Hm}, \ref{itm:Hf}
  and \ref{itm:Hlambda} hold, and set
\[
  m_0:=\max_{k\in\{1,\ldots,K\}}
  \max_{\alpha\in[\alpha_{2k-1},\alpha_{2k}]}m(\alpha).
\]
  Then $m_0>0$ and there exist $\varepsilon_1>0$ and
  $\eta_0\in(0,t_*)$ such that
\begin{equation}\label{eq:Hlambda-uniform-small}
  f(x,s)\geq
  (\lambda_1m_0+\varepsilon_1)s^{p-1}
  >
  \lambda_1m_0s^{p-1}
  \quad
  \hbox{for a.e. }x\in\Omega
  \hbox{ and every }s\in(0,\eta_0].
\end{equation}
  Consequently, for every $\theta_0\in(0,\eta_0)$ and every
  $\alpha\in\bigcup_{k=1}^K[\alpha_{2k-1},\alpha_{2k}]$,
  the function $\theta_0e_1$ is a subsolution of
\begin{equation}\label{eq:frozen-problem-general}
  -m(\alpha)\Delta_p u=f(x,u),\qquad u\in W^{1,p}_0(\Omega),
\end{equation}
  that is,
\[
  m(\alpha)\lambda_1(\theta_0e_1)^{p-1}\leq f(x,\theta_0e_1)
  \quad\hbox{for a.e. }x\in\Omega.
\]
\end{lemma}

\begin{proof}
Assumption \ref{itm:Hm} gives $m_0>0$. By \ref{itm:Hlambda}, there are
$\varepsilon_1>0$ and $\eta_0\in(0,t_*)$ such that
\[
 \essinf_{x\in\Omega}\frac{f(x,s)}{s^{p-1}}
 >\lambda_1m_0+\varepsilon_1
 \qquad\hbox{for every }s\in(0,\eta_0].
\]
The null set argument recalled at the beginning of this
section gives \eqref{eq:Hlambda-uniform-small}. For $0<\theta_0<\eta_0$, the
inequalities $0<e_1\leq1$ and $m(\alpha)\leq m_0$ then yield
\[
 -m(\alpha)\Delta_p(\theta_0e_1)
 =m(\alpha)\lambda_1(\theta_0e_1)^{p-1}
 \leq f(x,\theta_0e_1)
 \quad\hbox{a.e. in }\Omega.
\]
\end{proof}

\begin{lemma}\label{lem:large-barrier-near-ends}
  Assume that \ref{itm:Hm}, \ref{itm:Hf} 
  and \ref{itm:Hlambda} hold.
  Let $\theta_*\in(0,t_*)$. There exists $c_0>0$ such that
\begin{equation}\label{eq:cstar-bound}
  f(x,s)\geq \lambda_1c_0s^{p-1}
  \quad\hbox{for a.e. }x\in\Omega\hbox{ and every }s\in(0,\theta_*].
\end{equation}
  Consequently, if $m(\alpha)<c_0$ and $\theta\in(0,\theta_*]$, 
  then $\theta e_1$ is a subsolution 
  of \eqref{eq:frozen-problem-general}.
\end{lemma}

\begin{proof}
Use $m_0$ and $\eta_0$ from Lemma~\ref{lem:small-barrier}. If
$\theta_*\leq\eta_0$, take $c_0=m_0$. Otherwise,
\ref{itm:Hf}(ii) applied to $J=[\eta_0,\theta_*]$ and the null set
argument give $\rho_*>0$ such that
\[
 f(x,s)\geq\rho_*\geq\frac{\rho_*}{\theta_*^{p-1}}s^{p-1}
 \quad\hbox{for a.e. }x\in\Omega\hbox{ and every }s\in[\eta_0,\theta_*].
\]
Together with \eqref{eq:Hlambda-uniform-small}, this proves
\eqref{eq:cstar-bound} with
\[
 c_0=\min\left\{m_0,\frac{\rho_*}{\lambda_1\theta_*^{p-1}}\right\}>0.
\]
If $m(\alpha)<c_0$ and $0<\theta\leq\theta_*$, substituting
$s=\theta e_1(x)$ in \eqref{eq:cstar-bound} gives the subsolution inequality.
\end{proof}

\begin{lemma}[A continuous path of subsolutions]
\label{lem:plateau-family}
Assume \ref{itm:Hmu}, \ref{itm:Hm}, \ref{itm:Hf},
\ref{itm:Hlambda} and \ref{itm:He}. There exist
$0<\theta_0<\theta_*<t_*$, $w\in C^1_0(\overline\Omega)$,
and a continuous map
$[0,1]\ni r\mapsto L_r\in C^1_0(\overline\Omega)$ such that
\begin{equation}\label{eq:plateau-family-ends}
 L_0=\theta_0e_1,\qquad L_1=w,\qquad
 0<L_r\leq\theta_*\text{ in }\Omega,\qquad
 \mathcal H(w)>\alpha_{2K}.
\end{equation}
The function $L_0$ is a subsolution of \eqref{eq:frozen-problem-general}
throughout the positivity intervals of $m$.
There exists $c_*>0$ such that, for every $r\in[0,1]$,
$L_r$ is a subsolution of \eqref{eq:frozen-problem-general} whenever
$0<m(\alpha)<c_*$.
\end{lemma}
\begin{proof}
By \ref{itm:He} and the continuity of
$[0,t_*]\ni\theta\mapsto\mathcal H(\theta)$, choose $\theta_*\in(0,t_*)$ with
$\mathcal H(\theta_*)>\alpha_{2K}$. Fix an integer
$n>\max\{1,1/(p-1)\}$ and put
\[
 \zeta(s)=\begin{cases}
 1,&0\leq s\leq1,\\
 (2-s)^n,&1<s<2,\\
 0,&s\geq2,
 \end{cases}
 \qquad A=\int_0^2\zeta(s)\,ds=\frac{n+2}{n+1}.
\]
For $0<\delta<1/2$, let
\[
 h_\delta(s)=\frac1A\int_0^{s/\delta}\zeta(t)\,dt
 \quad(s\geq0).
\]
This is a concave, nondecreasing $C^1$ function, linear with slope
$(A\delta)^{-1}$ on $[0,\delta]$, and equal to $1$ on $[2\delta,\infty)$.
In particular $h_\delta(0)=0$, $h_\delta(1)=1$, and
$s\leq h_\delta(s)\leq1$ for $s\in[0,1]$. As $e_1>0$ in $\Omega$,
$h_\delta(e_1(x))\to1$ for every $x\in\Omega$ as $\delta\to0^+$.
Dominated convergence, applied to the weighted integrals, yields
\[
 \mathcal H(\theta_*h_\delta(e_1))\longrightarrow
 \mathcal H(\theta_*)>\alpha_{2K}.
\]
Fix $\delta$ so small that $w=\theta_*h_\delta(e_1)$ satisfies
$\mathcal H(w)>\alpha_{2K}$. Only convergence of these integrals is used. Uniform convergence up to the boundary is neither asserted nor needed.

For $r\in[0,1]$ define
\[
 \psi_r(s)=(1-r)s+rh_\delta(s),\qquad
 v_r=\psi_r(e_1),\qquad a_r(s)=(\psi_r'(s))^{p-1}.
\]
Then $e_1\leq v_r\leq1$ and $r\mapsto v_r$ is continuous in
$C^1_0(\overline\Omega)$. We verify a bound on $a_r'$ uniform in $r$. This also addresses the singular case $1<p<2$. Set
$B=\max\{1,(A\delta)^{-1}\}$. Clearly $0\leq a_r\leq B^{p-1}$.
On $\delta<s<2\delta$,
\[
 h_\delta'(s)=\frac{(2-s/\delta)^n}{A\delta},\qquad
 h_\delta''(s)=-\frac{n(2-s/\delta)^{n-1}}{A\delta^2}.
\]
For $p\geq2$ it follows that
$|a_r'|\leq(p-1)B^{p-2}n/(A\delta^2)$. If $1<p<2$ and $r>0$,
the negative exponent $p-2$ and $1-r+rh_\delta'\geq rh_\delta'$ give
\begin{align*}
 |a_r'(s)|
 &=(p-1)r(1-r+rh_\delta'(s))^{p-2}(-h_\delta''(s))\\
 &\leq\frac{(p-1)n r^{p-1}}{A^{p-1}\delta^p}
       (2-s/\delta)^{n(p-1)-1}
 \leq\frac{(p-1)n}{A^{p-1}\delta^p}.
\end{align*}
For $r=0$ the derivative vanishes, as it does outside the transition
interval. Since $n(p-1)>1$, $a_r$ is Lipschitz across both joining
points, also for $r=1$. Thus $a_r'\leq0$ and $|a_r'|\leq D$, uniformly,
where
\[
 D=\begin{cases}
 (p-1)B^{p-2}n/(A\delta^2),&p\geq2,\\
 (p-1)n/(A^{p-1}\delta^p),&1<p<2.
 \end{cases}
\]
Testing the weak eigenfunction equation with $a_r(e_1)\varphi$ gives
\begin{equation}\label{eq:plateau-chain}
 -\Delta_pv_r=\lambda_1a_r(e_1)e_1^{p-1}
              -a_r'(e_1)|\nabla e_1|^p
 \quad\text{in the weak sense}.
\end{equation}
Indeed, $\nabla v_r=\psi_r'(e_1)\nabla e_1$, and the product test is
admissible since $a_r$ is Lipschitz and $e_1\in C^1(\overline\Omega)$.
The chain rule is understood a.e. The values of $a_r'$ at joining points
are immaterial since $\nabla e_1=0$ a.e. on each corresponding level set.
No second derivatives of $e_1$ are required.
On $\{e_1\leq\delta\}$, $v_r$ is a constant multiple of $e_1$, so the
right-hand side is $\lambda_1v_r^{p-1}$. On $\{e_1>\delta\}$,
$v_r\geq e_1>\delta$. Consequently,
\begin{equation}\label{eq:plateau-differential-bound}
 0\leq-\Delta_pv_r\leq C_\delta v_r^{p-1},\qquad
 C_\delta=\max\left\{\lambda_1,
 \frac{\lambda_1B^{p-1}+D\|\nabla e_1\|_\infty^p}{\delta^{p-1}}\right\}.
\end{equation}
Choose $0<\theta_0<\min\{\theta_*,\eta_0\}$, with $\eta_0$ from
Lemma~\ref{lem:small-barrier}, and set
\[
 L_r=((1-r)\theta_0+r\theta_*)v_r.
\]
By homogeneity, \eqref{eq:plateau-differential-bound} holds also with
$L_r$ in place of $v_r$. Lemma~\ref{lem:large-barrier-near-ends} gives
$f(x,s)\geq\lambda_1c_0s^{p-1}$ for $0<s\leq\theta_*$. Hence the choice
\begin{equation}\label{eq:flat-threshold}
 c_*=\frac{\lambda_1c_0}{C_\delta}>0
\end{equation}
ensures $-m(\alpha)\Delta_pL_r\leq f(x,L_r)$ whenever
$0<m(\alpha)<c_*$. The remaining assertions follow from
the construction and Lemma~\ref{lem:small-barrier}.
\end{proof}

\section{Auxiliary problems and endpoint estimates}
\label{SecAux}
Throughout this section the assumptions of Theorem~\ref{thm:main} are
in force. We first record the properties of the functional
$\mathcal H$ defined in \eqref{eq:H-functional}.

\begin{lemma}\label{lem:H-properties}
Under \ref{itm:Hmu}, the following properties hold.
  \begin{enumerate}
    \item[(a)] $\mathcal H$ is continuous on $C(\overline\Omega)$,
in particular on $C^1_0(\overline\Omega)$.
    \item[(b)] 
$\mathcal H$ preserves order on nonnegative functions, in the sense that
      $0\leq v\leq u$ a.e. in $\Omega$ implies
      $\mathcal H(v)\leq\mathcal H(u)$.
  \end{enumerate}
\end{lemma}
\begin{proof}
  (a) Let $u_n\to u$ in $C(\overline\Omega)$ and set
  $c:=\sup_n\Vert u_n\Vert_\infty+\Vert u\Vert_\infty<+\infty$. Since
  $q_i\geq1$, the function $t\mapsto|t|^{q_i}$ is Lipschitz on $[-c,c]$,
  say with constant $L_i$, whence
\[
  |\mathcal H(u_n)-\mathcal H(u)|
  \leq\sum_{i=1}^{M}L_i\Vert\mu_i\Vert_\infty|\Omega|\,
  \Vert u_n-u\Vert_\infty\longrightarrow0 .
\]
  (b) Immediate from $\mu_i\geq0$ a.e. in $\Omega$ and the monotonicity
  of $t\mapsto t^{q_i}$ on $[0,+\infty)$.
\end{proof}

Fix $\theta_0$, $\theta_*$, $w$, $c_*$ and the family
$L_r$ given by Lemma~\ref{lem:plateau-family}. In particular,
\begin{equation}\label{eq:theta-star-large-H}
 \mathcal H(w)>\alpha_{2K}.
\end{equation}
We next arrange these barriers continuously along each parameter interval.

\begin{lemma}\label{lem:choice-intervals-barriers}
For each $k\in\{1,\ldots,K\}$ there exist
\begin{equation}\label{eq:beta-gamma-choice}
 \beta_k\in(\alpha_{2k-1},\tau_k),\qquad
 \gamma_k\in(\tau_k,\alpha_{2k})
\end{equation}
and continuous maps $\ell_k^-:[\beta_k,\tau_k]\to C^1_0(\overline\Omega)$,
$\ell_k^+:[\tau_k,\gamma_k]\to C^1_0(\overline\Omega)$ such that
$0<\ell_k^\pm(\alpha)\leq\theta_*<t_*$ in $\Omega$ and
\begin{equation}\label{eq:barrier-subsolutions}
 -m(\alpha)\Delta_p\ell_k^\pm(\alpha)
 \leq f(x,\ell_k^\pm(\alpha))\quad\text{in the weak sense}.
\end{equation}
Moreover,
\begin{equation}\label{eq:barrier-end-values}
 \ell_k^-(\beta_k)=\ell_k^+(\gamma_k)=w,\qquad
 \ell_k^-(\tau_k)=\ell_k^+(\tau_k)=\theta_0e_1.
\end{equation}
\end{lemma}
\begin{proof}
Continuity of $m$ and its zero endpoint values allow us to choose
\[
 \alpha_{2k-1}<\beta_k<\rho_k^-<\tau_k
 <\rho_k^+<\gamma_k<\alpha_{2k}
\]
such that
\begin{equation}\label{eq:m-small-near-ends}
 m(\alpha)<c_*\quad\text{for }
 \alpha\in[\alpha_{2k-1},\rho_k^-]\cup[\rho_k^+,\alpha_{2k}].
\end{equation}
Define
\[
 r_k^-(\alpha)=\begin{cases}
 (\rho_k^--\alpha)/(\rho_k^--\beta_k),&\beta_k\leq\alpha\leq\rho_k^-,\\
 0,&\rho_k^-\leq\alpha\leq\tau_k,
 \end{cases}
\]
\[
 r_k^+(\alpha)=\begin{cases}
 0,&\tau_k\leq\alpha\leq\rho_k^+,\\
 (\alpha-\rho_k^+)/(\gamma_k-\rho_k^+),&\rho_k^+\leq\alpha\leq\gamma_k,
 \end{cases}
\]
and put $\ell_k^\pm(\alpha)=L_{r_k^\pm(\alpha)}$. On the outer strips,
\eqref{eq:m-small-near-ends} and Lemma~\ref{lem:plateau-family} give the
subsolution inequality for all values of $r$. On the remaining parts,
$r_k^\pm=0$, so $\ell_k^\pm=\theta_0e_1$, which is a subsolution
throughout the positivity intervals by Lemma~\ref{lem:small-barrier}.
Continuity and the endpoint values follow directly. For the illustration
below, write the amplitude factors as
\begin{equation}\label{eq:theta-minus}
 \theta_k^\pm(\alpha)=(1-r_k^\pm(\alpha))\theta_0
                         +r_k^\pm(\alpha)\theta_*.
\end{equation}
Thus $\ell_k^\pm(\alpha)=\theta_k^\pm(\alpha)
\psi_{r_k^\pm(\alpha)}(e_1)$. Both the amplitude and the spatial profile
vary on the outer strips.
\end{proof}

\begin{figure}[htb]
\centering
\begin{tikzpicture}[x=0.75cm,y=0.75cm,>=Stealth,font=\small]

\def\xaxisleft{0.0}
\def\xa{1.2}
\def\xb{2.7}
\def\xd{5.3}
\def\xe{8.3}
\def\xf{11.3}
\def\xh{13.9}
\def\xJ{15.4}

\def\Ytop{6.1}
\def\Ybot{0.0}

\def\cstar{1.8}
\def\thetazero{0.9}
\def\thetastar{2.2}

\draw[->] (\xaxisleft,\Ytop) -- (\xJ+0.8,\Ytop) node[right] {$\alpha$};
\draw[->] (\xaxisleft,\Ytop-0.4) -- (\xaxisleft,\Ytop+3.4) node[above] {$m(\alpha)$};

\draw[dashed] (\xaxisleft,\Ytop+\cstar) -- (\xJ,\Ytop+\cstar);
\node[left] at (\xaxisleft,\Ytop+\cstar) {$c_*$};

\draw[thick,smooth]
(\xa,\Ytop)
.. controls (2.0,\Ytop+0.35) and (3.0,\Ytop+0.55) ..
(\xd,\Ytop+1.68)
.. controls (6.2,\Ytop+2.35) and (7.2,\Ytop+2.85) ..
(\xe,\Ytop+3.00)
.. controls (9.2,\Ytop+2.90) and (10.0,\Ytop+2.45) ..
(10.8,\Ytop+1.95)
.. controls (11.0,\Ytop+1.82) and (11.1,\Ytop+1.73) ..
(\xf,\Ytop+1.66)
.. controls (12.2,\Ytop+1.10) and (13.3,\Ytop+0.55) ..
(\xJ,\Ytop);

\node at (9.0,\Ytop+3.2) {$m$};

\node at (2.8,\Ytop+1.45) {$m(\alpha)<c_*$};
\node at (13.6,\Ytop+1.45) {$m(\alpha)<c_*$};

\foreach \x in {\xa,\xb,\xd,\xe,\xf,\xh,\xJ}
{
  \draw (\x,\Ytop+0.10) -- (\x,\Ytop-0.10);
}

\node[below]     at (\xa,\Ytop-0.10) {$\alpha_{2k-1}$};
\node[below]     at (\xb+0.3,\Ytop) {$\beta_k$};
\node[below]     at (\xd+0.3,\Ytop+0.1) {$\rho_k^-$};
\node[below]     at (\xe+0.3,\Ytop) {$\tau_k$};
\node[below]     at (\xf-0.3,\Ytop+0.1) {$\rho_k^+$};
\node[below]     at (\xh-0.3,\Ytop) {$\gamma_k$};
\node[below]     at (\xJ,\Ytop-0.10) {$\alpha_{2k}$};

\fill (\xa,\Ytop) circle (1.4pt);
\fill (\xb,\Ytop+0.55) circle (1.4pt);
\fill (\xd,\Ytop+1.68) circle (1.4pt);
\fill (\xe,\Ytop+3.00) circle (1.4pt);
\fill (\xf,\Ytop+1.66) circle (1.4pt);
\fill (\xh,\Ytop+0.44) circle (1.4pt);
\fill (\xJ,\Ytop) circle (1.4pt);

\draw[->] (\xaxisleft,\Ybot) -- (\xJ+0.8,\Ybot) node[right] {$\alpha$};
\draw[->] (\xaxisleft,\Ybot-0.4) -- (\xaxisleft,\Ybot+3.1) node[above] {$\theta$};

\draw[dashed] (\xaxisleft,\Ybot+\thetazero) -- (\xJ,\Ybot+\thetazero);
\draw[dashed] (\xaxisleft,\Ybot+\thetastar) -- (\xJ,\Ybot+\thetastar);
\node[left] at (\xaxisleft,\Ybot+\thetazero) {$\theta_0$};
\node[left] at (\xaxisleft,\Ybot+\thetastar) {$\theta_*$};

\draw[very thick,red]
(\xb,\Ybot+\thetastar)
-- (\xd,\Ybot+\thetazero)
-- (\xe,\Ybot+\thetazero);

\node[red] at (4.9,\Ybot+2.55) {$\theta_k^-$};

\draw[very thick,blue]
(\xe,\Ybot+\thetazero)
-- (\xf,\Ybot+\thetazero)
-- (\xh,\Ybot+\thetastar);

\node[blue] at (11.8,\Ybot+2.55) {$\theta_k^+$};

\foreach \x in {\xa,\xb,\xd,\xe,\xf,\xh,\xJ}
{
  \draw (\x,\Ybot+0.10) -- (\x,\Ybot-0.10);
}

\node[below]     at (\xa,\Ybot-0.10) {$\alpha_{2k-1}$};
\node[below]     at (\xb,\Ybot-0.10) {$\beta_k$};
\node[below]     at (\xd,\Ybot-0.10) {$\rho_k^-$};
\node[below]     at (\xe,\Ybot-0.10) {$\tau_k$};
\node[below]     at (\xf,\Ybot-0.10) {$\rho_k^+$};
\node[below]     at (\xh,\Ybot-0.10) {$\gamma_k$};
\node[below]     at (\xJ,\Ybot-0.10) {$\alpha_{2k}$};

\draw[dashed] (\xb,\Ybot) -- (\xb,\Ytop+0.8);
\draw[dashed] (\xd,\Ybot) -- (\xd,\Ytop+1.7);
\draw[dashed] (\xe,\Ybot) -- (\xe,\Ytop+3);
\draw[dashed] (\xf,\Ybot) -- (\xf,\Ytop+1.7);
\draw[dashed] (\xh,\Ybot) -- (\xh,\Ytop+0.5);

\end{tikzpicture}
\caption{Construction of the parameter-dependent barriers in
Lemma~\ref{lem:choice-intervals-barriers}. The upper panel shows a schematic
positivity interval of $m$. The lower panel shows the amplitude factors
$\theta_k^-$ and $\theta_k^+$. The spatial profile also varies through
$\psi_{r_k^\pm}(e_1)$. No unimodality assumption on $m$ is used.}
\label{fig:lemma31-geometry}
\end{figure}
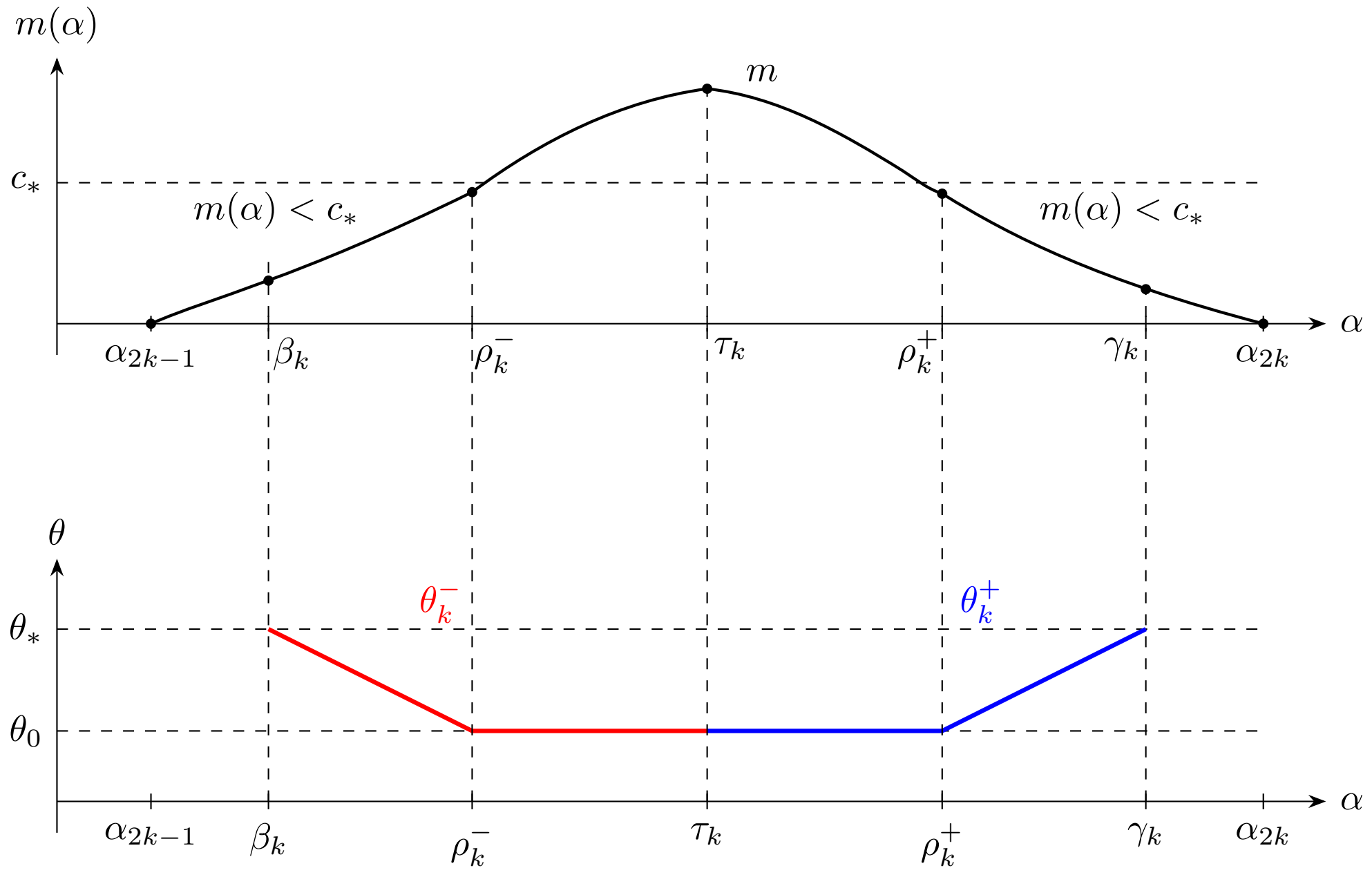

Figure~\ref{fig:lemma31-geometry} illustrates the construction.

  From now on we fix, for every $k\in\{1,\ldots,K\}$, a choice of
  $\beta_k$, $\gamma_k$ and of the barrier functions $\ell_k^\pm$
  provided by Lemma \ref{lem:choice-intervals-barriers}.
  Let $I$ be one of the two intervals
\[
  I_k^-:=[\beta_k,\tau_k],
  \qquad
  I_k^+:=[\tau_k,\gamma_k],
\]
  and let $\ell_I(\alpha)$ denote 
  the corresponding lower barrier, namely
\[
  \ell_{I_k^-}(\alpha)=\ell_k^-(\alpha),
  \qquad
  \ell_{I_k^+}(\alpha)=\ell_k^+(\alpha).
\]
  For $\alpha\in I$ define the truncation
\begin{equation}\label{eq:truncation}
  \widehat f_I(\alpha,x,s)=
  \left\{
  \begin{array}{ll}
    f(x,\ell_I(\alpha)(x)), & s<\ell_I(\alpha)(x),\\
    f(x,s), & \ell_I(\alpha)(x)\leq s\leq t_*,\\
    0, & s>t_*.
  \end{array}
  \right.
\end{equation}

\begin{remark}\label{rem:hatf-caratheodory}
For every $(\alpha,s)\in I\times\R$ the function
$\widehat f_I(\alpha,\cdot,s)$ is measurable, and for almost every
$x\in\Omega$ the function $\widehat f_I(\cdot,x,\cdot)$ is continuous on
$I\times\R$. Indeed, $\alpha\mapsto\ell_I(\alpha)$ is continuous with
values in $C^1_0(\overline\Omega)$, and the pieces in
\eqref{eq:truncation} agree at both joining points, at $s=t_*$ because
of \ref{itm:Hf}(i). Moreover,
\eqref{eq:f-upper-bound-uniform} gives
\begin{equation}\label{eq:rhs-bound}
 0\leq\widehat f_I(\alpha,x,s)\leq T
 \quad\hbox{for a.e. }x\in\Omega\hbox{ and every }(\alpha,s)\in I\times\R.
\end{equation}
\end{remark}

\begin{lemma}\label{lem:solution-operator}
  Fix $k\in\{1,\ldots,K\}$ and let $I\in\{I_k^-,I_k^+\}$.
  For every $(\alpha,v)\in I\times C^1_0(\overline\Omega)$,
  the problem
\begin{equation}\label{eq:T-problem}
  \left\{
  \begin{array}{ll}
    -m(\alpha)\Delta_p z=\widehat f_I(\alpha,x,v(x)) 
      &\hbox{in }\Omega,\\
    z=0 &\hbox{on }\partial\Omega
  \end{array}\right.
\end{equation}
  has a unique weak solution $z$, 
  and $z\geq0$ a.e. in $\Omega$.
  We denote it by
\[
  \mathcal T_I(\alpha,v):=z.
\]
  Moreover,
  $\mathcal T_I:I\times C^1_0(\overline\Omega)\to
  C^1_0(\overline\Omega)$
  is continuous, and 
  $\mathcal T_I\big(I\times C^1_0(\overline\Omega)\big)$
  is a relatively compact subset 
  of $C^1_0(\overline\Omega)$.
\end{lemma}
\begin{proof}
The interval $I$ is a compact subset of a positivity interval of $m$, so
\begin{equation}\label{eq:mI-positive}
 m_I:=\min_{\alpha\in I}m(\alpha)>0.
\end{equation}
For fixed $(\alpha,v)$ the right-hand side of
\eqref{eq:T-problem} is a prescribed function in $L^\infty(\Omega)$,
bounded between $0$ and $T$ by \eqref{eq:rhs-bound}. The theorem of
Minty and Browder therefore gives a unique weak solution
$z$. Testing with $-z^-$ gives
\[
 m(\alpha)\int_\Omega|\nabla z^-|^p\,dx
 =-\int_\Omega\widehat f_I(\alpha,x,v)z^-\,dx\leq0,
\]
and hence $z\geq0$.

Next, let $w_0\in W^{1,p}_0(\Omega)\cap L^\infty(\Omega)$
be the solution of $-\Delta_pw_0=1$ in $\Omega$. Since
\[
 -\Delta_pz=\frac{\widehat f_I(\alpha,x,v)}{m(\alpha)}\leq\frac{T}{m_I}
 =-\Delta_p\Bigl(\Bigl(\frac{T}{m_I}\Bigr)^{1/(p-1)}w_0\Bigr)
 \quad\hbox{in }\Omega,
\]
the weak comparison principle gives $0\leq z\leq(T/m_I)^{1/(p-1)}w_0$,
a bound in $L^\infty(\Omega)$ independent of $(\alpha,v)$. The global
$C^{1,\nu}$ regularity estimate for bounded data \cite{Lieberman}
therefore gives
$\nu\in(0,1)$ and $C_I>0$, independent of $(\alpha,v)$, such that
\begin{equation}\label{eq:regularity-bound}
 \Vert\mathcal T_I(\alpha,v)\Vert_{C^{1,\nu}(\overline\Omega)}\leq C_I.
\end{equation}
The compact embedding into $C^1(\overline\Omega)$ proves relative
compactness of the image.

To prove continuity, let $(\alpha_n,v_n)\to(\alpha,v)$ in
$I\times C^1_0(\overline\Omega)$ and set $z_n=\mathcal T_I(\alpha_n,v_n)$.
Every subsequence has a further subsequence converging in
$C^1_0(\overline\Omega)$ to some $z$. By
Remark~\ref{rem:hatf-caratheodory} and dominated convergence,
\[
 \widehat f_I(\alpha_n,x,v_n(x))
 \longrightarrow\widehat f_I(\alpha,x,v(x))
 \quad\hbox{in }L^{p'}(\Omega).
\]
Passing to the limit in the weak formulation shows that
$z=\mathcal T_I(\alpha,v)$ by uniqueness for the prescribed
right-hand side.
Every convergent subsequence has this same limit, so the full sequence
converges in $C^1_0(\overline\Omega)$.
\end{proof}

Uniqueness in the above lemma concerns the problem with the given input $v$.
It does not assert uniqueness of fixed points of
$\mathcal T_I(\alpha,\cdot)$.

\begin{lemma}\label{lem:order-bounds-fixed-points}
  Fix $k\in\{1,\ldots,K\}$, 
  let $I\in\{I_k^-,I_k^+\}$ and let
  $(\alpha,u)\in I\times C^1_0(\overline\Omega)$ 
  be such that
\[
  u=\mathcal T_I(\alpha,u).
\]
  Then
\begin{equation}\label{eq:order-bounds}
  \ell_I(\alpha)\leq u\leq t_*\quad\hbox{a.e. in }\Omega.
\end{equation}
  Consequently, $u$ solves the untruncated frozen problem
\begin{equation}\label{eq:untruncated-frozen}
  -m(\alpha)\Delta_pu=f(x,u),
  \qquad u\in W^{1,p}_0(\Omega).
\end{equation}
\end{lemma}
\begin{proof}
Write $\ell=\ell_I(\alpha)$. Since $u=\mathcal T_I(\alpha,u)$,
\begin{equation}\label{eq:u-truncated-weak}
 m(\alpha)\int_\Omega|\nabla u|^{p-2}\nabla u\cdot\nabla\varphi\,dx
 =\int_\Omega\widehat f_I(\alpha,x,u)\varphi\,dx
 \quad\hbox{for every }\varphi\in W^{1,p}_0(\Omega).
\end{equation}
Take $\varphi=(\ell-u)^+$. On the set $\{u<\ell\}$, where
$\varphi$ may be nonzero, the truncated reaction equals $f(x,\ell)$.
Subtracting \eqref{eq:u-truncated-weak} from the weak form of the
subsolution inequality \eqref{eq:barrier-subsolutions} for $\ell$, both
tested with $\varphi$, yields
\[
 m(\alpha)\int_\Omega
 \bigl(|\nabla\ell|^{p-2}\nabla\ell-|\nabla u|^{p-2}\nabla u\bigr)
 \cdot\nabla(\ell-u)^+\,dx\leq0.
\]
The integrand is nonnegative by the monotonicity of
$\xi\mapsto|\xi|^{p-2}\xi$, and since $m(\alpha)>0$ it vanishes a.e. in
$\Omega$. By strict monotonicity, $\nabla\ell=\nabla u$ a.e. on
$\{u<\ell\}$, that is, $\nabla(\ell-u)^+=0$ a.e. in $\Omega$.
Poincar\'e's inequality gives
$(\ell-u)^+=0$, and hence $u\geq\ell$.

Next take $\varphi=(u-t_*)^+$, which belongs to $W^{1,p}_0(\Omega)$.
The truncated reaction vanishes on $\{u>t_*\}$, so
\[
 m(\alpha)\int_\Omega|\nabla(u-t_*)^+|^p\,dx=0.
\]
Poincar\'e's inequality now gives $u\leq t_*$. The truncation therefore
agrees with $f(x,u)$ a.e., proving \eqref{eq:untruncated-frozen}.
\end{proof}

\begin{lemma}\label{lem:midpoint-inequality}
Fix $k\in\{1,\ldots,K\}$. Every weak solution
$u\in W^{1,p}_0(\Omega)$ of
\begin{equation}\label{eq:midpoint-frozen}
 -m(\tau_k)\Delta_pu=f(x,u)\quad\hbox{in }\Omega,
 \qquad u=0\quad\hbox{on }\partial\Omega,
\end{equation}
with $0\leq u\leq t_*$ a.e. satisfies
\begin{equation}\label{eq:midpoint-H-less}
 \mathcal H(u)<\tau_k.
\end{equation}
\end{lemma}

\begin{proof}
Testing with $u$ and using \eqref{eq:tf-upper-bound-uniform} gives
\[
 m(\tau_k)\Vert u\Vert^p
 =\int_\Omega f(x,u)u\,dx\leq|\Omega|E.
\]
Since $u^{q_i}\leq t_*^{q_i-1}u$ for $0\leq u\leq t_*$,
H\"older's inequality and \eqref{eq:lambda1} imply
\begin{align}\label{eq:H-midpoint-estimate}
 \mathcal H(u)
 &\leq\sum_{i=1}^M t_*^{q_i-1}\Vert\mu_i\Vert_{p'}\Vert u\Vert_p
 \leq\lambda_1^{-1/p}
       \sum_{i=1}^M t_*^{q_i-1}\Vert\mu_i\Vert_{p'}\Vert u\Vert
 \notag\\
 &\leq\left(\frac{|\Omega|E}{\lambda_1m(\tau_k)}\right)^{1/p}
       \sum_{i=1}^M t_*^{q_i-1}\Vert\mu_i\Vert_{p'}
 <\tau_k.
\end{align}
The last inequality is \ref{itm:Htau} in equivalent form.
\end{proof}

\section{Proofs of the main results}
\label{SecMain}

\begin{proof}[Proof of Theorem~\ref{thm:main}]
Set $X=C^1_0(\overline\Omega)$. The uniform bounds in
\eqref{eq:regularity-bound} and the finiteness of the family of intervals
$I_k^\pm$ give $R>0$ such that each $\mathcal T_{I_k^\pm}$ maps
$I_k^\pm\times B_R$ into
$B_R=\{v\in X:\Vert v\Vert_{C^1(\overline\Omega)}\leq R\}$.
By Lemmas~\ref{lem:H-properties} and \ref{lem:solution-operator},
$\Phi=\mathcal H$ is continuous and these solution operators are
continuous and compact. Thus the functional and operator assumptions of
Lemma~\ref{lem:abstract-constraint-fp} hold.

Fix $k$. If $u=\mathcal T_{I_k^-}(\beta_k,u)$ or
$u=\mathcal T_{I_k^+}(\gamma_k,u)$, then
Lemma~\ref{lem:order-bounds-fixed-points} and
\eqref{eq:barrier-end-values} give $u\geq w$. Since
$\mathcal H$ preserves order on nonnegative functions,
\begin{equation}\label{eq:outer-endpoint-bound}
 \mathcal H(u)\geq\mathcal H(w)
 >\alpha_{2K}\geq\alpha_{2k}>\gamma_k>\beta_k.
\end{equation}
At the common endpoint $\tau_k$, every fixed point of either operator
solves \eqref{eq:midpoint-frozen} and lies between $0$ and $t_*$ by
Lemma~\ref{lem:order-bounds-fixed-points}. Lemma~\ref{lem:midpoint-inequality}
therefore gives
\[
 \mathcal H(u)<\tau_k.
\]
These estimates hold for every fixed point at the
respective endpoints, as required in
Lemma~\ref{lem:abstract-constraint-fp}.

Apply Lemma~\ref{lem:abstract-constraint-fp}(a) on
$[\beta_k,\tau_k]$ and part (b) on $[\tau_k,\gamma_k]$.
We obtain pairs $(a_{k,1},u_{k,1})$ and $(a_{k,2},u_{k,2})$ with
\[
 u_{k,1}=\mathcal T_{I_k^-}(a_{k,1},u_{k,1}),\qquad
 u_{k,2}=\mathcal T_{I_k^+}(a_{k,2},u_{k,2}),\qquad
 a_{k,j}=\mathcal H(u_{k,j}),
\]
and
\[
 \beta_k<a_{k,1}<\tau_k<a_{k,2}<\gamma_k.
\]
By Lemma~\ref{lem:order-bounds-fixed-points}, each $u_{k,j}$ solves
\[
 -m(a_{k,j})\Delta_pu_{k,j}=f(x,u_{k,j}),\qquad
 0<u_{k,j}\leq t_*\quad\hbox{in }\Omega.
\]
Substituting $a_{k,j}=\mathcal H(u_{k,j})$ proves that these are positive
solutions of \eqref{Problem1}. They belong to $C^1_0(\overline\Omega)$
by construction. Since $\alpha_{2k-1}<\beta_k$ and
$\gamma_k<\alpha_{2k}$ by \eqref{eq:beta-gamma-choice}, this gives
\eqref{eq:main-order}. Distinct nonlocal values separate the two
solutions in each pair, and the disjoint positivity intervals of $m$
separate the pairs. Hence all $2K$ solutions are distinct.
\end{proof}

\begin{proof}[Proof of Theorem~\ref{thm:perturbed-functional}]
Fix the function $w$ given by Lemma~\ref{lem:plateau-family}
and the intervals $I_k^-=[\beta_k,\tau_k]$ and
$I_k^+=[\tau_k,\gamma_k]$, where $\beta_k$ and $\gamma_k$
are chosen as in 
Lemma~\ref{lem:choice-intervals-barriers}. 
Define
\begin{equation}\label{eq:perturb-Rk}
 R_k=\left(\frac{|\Omega|E}{\lambda_1m(\tau_k)}\right)^{1/p}
       \sum_{i=1}^M t_*^{q_i-1}\|\mu_i\|_{p'},
\end{equation}
\begin{equation}\label{eq:perturb-margin}
 d=\min\left\{\mathcal H(w)-\alpha_{2K},
              \min_{1\leq k\leq K}(\tau_k-R_k)\right\}>0,
 \qquad 0<\varepsilon_0<d.
\end{equation}
Positivity follows from \eqref{eq:theta-star-large-H} and
\ref{itm:Htau}. Any choice of $\varepsilon_0\in(0,d)$
gives a sufficient perturbation tolerance.

Fix a continuous functional
$\widetilde{\mathcal H}:X\to[0,\infty)$ satisfying
\eqref{eq:perturb-smallness}.

Use the barriers $\ell_k^\pm$ and operators
$\mathcal T_{I_k^\pm}$ constructed in Section~\ref{SecAux}. Their construction uses
$f$, $m$ and the fixed barriers, and does not require the functional
appearing in the scalar constraint to preserve order. Choose the same
ball $B_R\subset X$ as in the proof of Theorem~\ref{thm:main}.
Lemma~\ref{lem:solution-operator} gives continuous compact maps
$I_k^\pm\times B_R\to B_R$. We apply
Lemma~\ref{lem:abstract-constraint-fp} with
$\Phi=\widetilde{\mathcal H}$.

Every fixed point lies in $\mathcal D$ by
Lemma~\ref{lem:order-bounds-fixed-points}. At either outer endpoint,
\eqref{eq:barrier-end-values} gives $u\geq w$. The reference functional
$\mathcal H$ preserves order, so
\begin{align*}
 \widetilde{\mathcal H}(u)
 &\geq\mathcal H(u)-\varepsilon_0
 \geq\mathcal H(w)-\varepsilon_0\\
 &\geq\alpha_{2K}+d-\varepsilon_0
 >\alpha_{2K}>\gamma_k>\beta_k.
\end{align*}
At the common endpoint $\tau_k$, the fixed point solves the untruncated
frozen equation. The estimate \eqref{eq:H-midpoint-estimate} gives
$\mathcal H(u)\leq R_k$, and hence
\[
 \widetilde{\mathcal H}(u)\leq R_k+\varepsilon_0
 \leq\tau_k-d+\varepsilon_0<\tau_k.
\]
Both estimates hold for all fixed points at the corresponding endpoints.

Lemma~\ref{lem:abstract-constraint-fp}(a) on $I_k^-$ and part (b) on
$I_k^+$ therefore provide points
$(a_{k,j},u_{k,j})$ satisfying
\[
 \begin{gathered}
 u_{k,j}=\mathcal T_{I_k^\pm}(a_{k,j},u_{k,j}),\qquad
 a_{k,j}=\widetilde{\mathcal H}(u_{k,j}),\\
 \beta_k<a_{k,1}<\tau_k<a_{k,2}<\gamma_k.
 \end{gathered}
\]
By Lemma~\ref{lem:order-bounds-fixed-points}, substitution of this scalar
identity into the frozen equation gives \eqref{eq:perturbed-problem}.
Positivity follows from the lower barrier, and regularity follows from
the definition of $\mathcal T_I$. The disjoint intervals containing the
values of $\widetilde{\mathcal H}$ give \eqref{eq:perturbed-order} and
separate all $2K$ solutions.
\end{proof}

\section{Examples}\label{SecExample}

\begin{example}\label{ex:nonmonotone-nonlinearity}
Let $N=2$, $p=3$, $t_*>1$ and $\Omega=B(x_0,3)$. Choose
$x_2\in B(x_0,2)$, set $\Omega_1=B(x_0,1)$ and
$\Omega_2=B(x_2,1)$ (so that $\Omega_2\subset\Omega$), and define
\[
 \mu_1=\chi_{\Omega_1},\qquad \mu_2=2\chi_{\Omega_2}.
\]
These weights satisfy \ref{itm:Hmu}. The two regions may overlap, and
their union does not cover $\Omega$.

Fix any $q_1,q_2\geq1$. Since
$|\Omega_1|=|\Omega_2|=\pi$ and $t_*>1$,
\[
 M_*=\pi t_*^{q_1}+2\pi t_*^{q_2}>3\pi>\frac32.
\]
Thus \ref{itm:He} holds for the coefficient $m$ defined in
\eqref{eq:example-nonmonotone-m} below, for which $K=1$ and
$\alpha_{2K}=\alpha_2=3/2$.

  We consider the autonomous 
  Carath\'eodory reaction
  $f(x,t)=f(t)$, where
\[
  f(t)=\lambda_1\bigl((t_*-t)^+\bigr)^2
  \left(1+\frac34\sin\frac{4\pi t}{t_*}\right),
  \qquad t\in\R.
\]
  Since
\[
  1+\frac34\sin\frac{4\pi t}{t_*}\geq\frac14
  \qquad\hbox{for all }t\in\R,
\]
  we have $f(t)\geq0$ for all $t\in\R$. Moreover,
  $f(t_*)=0$. If $J\Subset(0,t_*)$, then
  $\inf\limits_{t\in J}f(t)>0$,
  and
\[
  \sup_{t\in[0,t_*]}f(t)
  \leq \frac74\lambda_1t_*^2=:T.
\]
  Thus \ref{itm:Hf} holds. Furthermore, since $p=3$,
\[
  \liminf_{t\to0^+}\frac{f(t)}{t^{p-1}}
  =
  \liminf_{t\to0^+}
  \lambda_1
  \frac{(t_*-t)^2}{t^2}
  \left(1+\frac34\sin\frac{4\pi t}{t_*}\right)
  =+\infty.
\]
Therefore \ref{itm:Hlambda} holds for every continuous
coefficient $m$ satisfying \ref{itm:Hm}, in particular for the one
defined in \eqref{eq:example-nonmonotone-m} below.

  Let us also observe that the quotient
  $t\mapsto f(t)/t^{p-1}$ is not monotone
  on $(0,t_*)$. Indeed, for $t\in(0,t_*)$ set
\[
  Q(t):=\frac{f(t)}{t^2}
  =
  \lambda_1
  \frac{(t_*-t)^2}{t^2}
  \left(1+\frac34\sin\frac{4\pi t}{t_*}\right).
\]
  Since $Q(t)\to+\infty$ as $t\to0^+$,
$Q$ is not nondecreasing.
  On the other hand, writing $Q=\lambda_1hg$ with
  $h(t)=\big(\frac{t_*-t}{t}\big)^2$ and
  $g(t)=1+\frac34\sin\frac{4\pi t}{t_*}$, we have
  $h(t_*/2)=g(t_*/2)=1$,
  $h'(t_*/2)=-8/t_*$ and 
  $g'(t_*/2)=3\pi/t_*$, whence
\[
  Q'\left(\frac{t_*}{2}\right)
  =\lambda_1\big(h'g+hg'\big)\left(\frac{t_*}{2}\right)
  =
  \frac{\lambda_1}{t_*}(3\pi-8)>0.
\]
Thus $Q$ is not nonincreasing either. This reaction violates the decrease
condition used for the multiplicity results in
\cite{CanGasLivSan,GasSan1,GasSan2,GasSanSic}.

  Let $\alpha_1=1<\alpha_2=\frac{3}{2}$
  and define
\begin{equation}\label{eq:example-nonmonotone-m}
  m(s)=As(s-1)\left(\frac32-s\right),
  \qquad s\geq0,
\end{equation}
  where
\begin{equation}\label{eq:example-nonmonotone-A}
  A>
  \frac{28672}{1875}\pi^3t_*^3
  \left(t_*^{q_1-1}+2t_*^{q_2-1}\right)^3.
\end{equation}
  Then $m>0$ in $(1,3/2)$ and $m(1)=m(3/2)=0$, hence \ref{itm:Hm}
  holds with $K=1$.

  It remains to verify \ref{itm:Htau}. 
  Since $K=1$, we choose
  $\tau_1=5/4$ and write $\tau:=\tau_1$ 
  for brevity.
  By \eqref{eq:example-nonmonotone-m} and
  \eqref{eq:example-nonmonotone-A},
\begin{equation}\label{eq:example-nonmonotone-taum}
  \tau^3m(\tau)
  =
  A\tau^4(\tau-1)\left(\frac32-\tau\right)
  =
  A\frac{625}{4096}
  >
  \frac73\pi^3t_*^3
  \left(t_*^{q_1-1}+2t_*^{q_2-1}\right)^3.
\end{equation}
  Moreover, for $t\in[0,t_*]$,
\[
  tf(t)
  \leq
  \frac74\lambda_1 t(t_*-t)^2.
\]
  Since
\[
  \max_{t\in[0,t_*]}t(t_*-t)^2
  =
  \frac{4t_*^3}{27},
\]
  we get
\[
  \max_{t\in[0,t_*]}[tf(t)]
  \leq
  \lambda_1\frac{7t_*^3}{27}.
\]
  Since $|\Omega|=9\pi$, $p'=3/2$, and
\[
  \Vert\mu_1\Vert_{3/2}
  =\Vert\chi_{\Omega_1}\Vert_{3/2}
  =\pi^{2/3},
  \qquad
  \Vert\mu_2\Vert_{3/2}
  =\Vert2\chi_{\Omega_2}\Vert_{3/2}
  =2\pi^{2/3},
\]
  we obtain
\begin{align*}
    &\frac{|\Omega|}{\lambda_1}
    \left(
      t_*^{q_1-1}\Vert\mu_1\Vert_{3/2}
      +
      t_*^{q_2-1}\Vert\mu_2\Vert_{3/2}
    \right)^3
    \max_{t\in[0,t_*]}[tf(t)]
    \\
    &\leq
    \frac{9\pi}{\lambda_1}
    \left(
      t_*^{q_1-1}\pi^{2/3}
      +
      2t_*^{q_2-1}\pi^{2/3}
    \right)^3
    \lambda_1\frac{7t_*^3}{27}
    =
    \frac73\pi^3t_*^3
    \left(t_*^{q_1-1}+2t_*^{q_2-1}\right)^3.
\end{align*}
  Combining this estimate 
  with \eqref{eq:example-nonmonotone-taum}, we obtain
  \ref{itm:Htau}.

All assumptions of Theorem~\ref{thm:main} are now verified. Problem
\eqref{Problem1}, with the weights, reaction and coefficient specified
above, therefore has two positive solutions
$u_1,u_2\in C^1_0(\overline\Omega)$ satisfying
\[
 1<\mathcal H(u_1)<\frac54<\mathcal H(u_2)<\frac32,
 \qquad
 \mathcal H(u)=\Vert u\Vert_{q_1,\Omega_1}^{q_1}
              +2\Vert u\Vert_{q_2,\Omega_2}^{q_2}.
\]
The normalized reaction is nonmonotone, so the cited multiplicity
theorems requiring a decreasing quotient do not apply to these data.

\end{example}

\begin{example}[Explicit data for a quadratic coefficient]
\label{ex:explicit-data}
Let $p=3$, $M=1$, $q_1=1$, $\mu_1=1$, and $t_*>0$, on any bounded
smooth domain as above. Then
 $\mathcal H(u)=\int_\Omega |u|\,dx$,
 $M_*=t_*|\Omega|$.
Set
$b=M_*/2$,
$a=b/2$,
$\tau=3b/4$,
and define
\[
 f(x,t)=(t_*-t)^2\quad(t\in\mathbb R),\qquad
 m(s)=A_0(s-a)(b-s)\quad(s\geq0).
\]
The reaction satisfies \ref{itm:Hf}, with $T=t_*^2$, and
$f(x,t)/t^2\to+\infty$ as $t\to0^+$, so \ref{itm:Hlambda} holds for
every finite $A_0>0$. The coefficient satisfies \ref{itm:Hm} with
$K=1$, $\alpha_1=a$ and $\alpha_2=b$. Since $b<M_*$,
\ref{itm:He} also holds. Moreover,
\[
 E=\max_{0\leq t\leq t_*}t(t_*-t)^2=\frac{4t_*^3}{27},
 \qquad \|\mu_1\|_{3/2}=|\Omega|^{2/3}.
\]
Consequently \ref{itm:Htau} holds whenever
\begin{equation}\label{eq:quadratic-coefficient-size}
 A_0>\frac{4|\Omega|^3t_*^3}
 {27\lambda_1\tau^3(\tau-a)(b-\tau)}.
\end{equation}
For every such $A_0$, Theorem~\ref{thm:main} gives two positive solutions
$u_1,u_2$ with
\[
 a<\mathcal H(u_1)<\tau<\mathcal H(u_2)<b.
\]
These data will also be used in Example~\ref{ex:nonorder-functional}.
\end{example}

\begin{example}[A perturbation depending on an additional spatial integral]
\label{ex:nonorder-functional}
This example illustrates Theorem~\ref{thm:perturbed-functional} through a
perturbation involving an additional spatial integral. We verify that
order preservation may fail and 
show that in general the effect of the perturbation cannot be reproduced
by changing $m$ while keeping $\mathcal H$ unchanged.

Use the data from Example~\ref{ex:explicit-data}, with $A_0$ satisfying
\eqref{eq:quadratic-coefficient-size}, and let $d>0$ be the margin
defined in \eqref{eq:perturb-margin} in the proof of
Theorem~\ref{thm:perturbed-functional}. Let $U\subseteq\Omega$ be any
measurable set with $|U|>0$, and put
\[
 J(u)=\int_U |u|\,dx,\qquad
 0<\eta<\min\{1/5,d/M_*\}.
\]
For $\omega>0$ define
\begin{equation}\label{eq:oscillating-functional}
 \widetilde{\mathcal H}_{\eta,\omega}(u)
 =\mathcal H(u)\bigl[1+\eta\sin(\omega J(u))\bigr],
 \qquad u\in X.
\end{equation}
This functional is continuous and nonnegative on $X$. For $u\in\mathcal D$,
\[
 |\widetilde{\mathcal H}_{\eta,\omega}(u)-\mathcal H(u)|
 \leq\eta\mathcal H(u)\leq\eta M_*<d.
\]
The proof of Theorem~\ref{thm:perturbed-functional} allows any
$\varepsilon_0\in(0,d)$. Choosing $\varepsilon_0\in(\eta M_*,d)$
therefore gives two positive solutions for every $\omega>0$.

We now choose $\omega$ so that the functional fails to preserve order.
Set $s_0=5b/8=5M_*/16$ and choose $v_0\in X$ such that
\[
 0<v_0<t_*\quad\text{in }\Omega,\qquad
 \mathcal H(v_0)=s_0.
\]
Such a choice is possible since $s_0<M_*$, by approximating the
constant $t_*$ from below with positive functions in $X$ and rescaling.
Since $|U|>0$ and $v_0>0$ in $\Omega$,
\[
 \rho=\frac{J(v_0)}{s_0}\in(0,1].
\]
Choose an integer $j\geq1$ such that $\eta(2j-1)\pi>1$, and set
\[
 \omega=\frac{(4j+1)\pi}{2\rho s_0},\qquad
 s_1=\frac{(2j-1)\pi}{\omega\rho}\in(0,s_0).
\]
Along the ordered family $u_s=(s/s_0)v_0$, $0<s<s_0$, we have
$\mathcal H(u_s)=s$, $J(u_s)=\rho s$ and $u_s\in\mathcal D$.
Consequently, for $G(s)=s[1+\eta\sin(\omega\rho s)]$,
\[
 \widetilde{\mathcal H}_{\eta,\omega}(u_s)=G(s),\qquad
 G'(s_1)=1-\eta(2j-1)\pi<0.
\]
There are therefore $s_-<s_+$ near $s_1$ with $u_{s_-}\leq u_{s_+}$ but
\[
 \widetilde{\mathcal H}_{\eta,\omega}(u_{s_-})
 >
 \widetilde{\mathcal H}_{\eta,\omega}(u_{s_+}).
\]
In particular, this functional cannot be a sum of integral powers with
nonnegative weights, since every functional of that form preserves order
on nonnegative functions.

To rule out absorption of the perturbation into $m$, assume in addition
that $|U|<|\Omega|$. Choose $h\in C^\infty_c(\Omega)$ such that
\[
 \int_\Omega h\,dx=0,\qquad \int_U h\,dx=1.
\]
Such an $h$ exists because the two integration functionals are linearly
independent on $C^\infty_c(\Omega)$. Otherwise the characteristic function
of $U$ would be constant a.e. in $\Omega$, contrary to $0<|U|<|\Omega|$.
For sufficiently small $\varepsilon>0$, the function
$v_\varepsilon=v_0+\varepsilon h$ belongs to $X$ and satisfies
$0<v_\varepsilon<t_*$ in $\Omega$. Moreover,
\[
 \mathcal H(v_\varepsilon)=\mathcal H(v_0)=s_0,\qquad
 J(v_\varepsilon)=J(v_0)+\varepsilon.
\]
Since $\omega J(v_0)=(4j+1)\pi/2$, choosing also
$0<\omega\varepsilon<\pi/2$ gives
\[
 \widetilde{\mathcal H}_{\eta,\omega}(v_\varepsilon)
 =s_0[1+\eta\cos(\omega\varepsilon)]
 <s_0(1+\eta)=\widetilde{\mathcal H}_{\eta,\omega}(v_0).
\]
As $\eta<1/5$, both values belong to $(a,\tau)$, since
\[
 a<s_0<\widetilde{\mathcal H}_{\eta,\omega}(v_\varepsilon)
 <\widetilde{\mathcal H}_{\eta,\omega}(v_0)
 =s_0(1+\eta)<\tau=(a+b)/2.
\]
The coefficient $m(s)=A_0(s-a)(b-s)$ is strictly increasing on
$(a,\tau)$. Hence both coefficient values are positive and
\[
 m\bigl(\widetilde{\mathcal H}_{\eta,\omega}(v_\varepsilon)\bigr)
 <
 m\bigl(\widetilde{\mathcal H}_{\eta,\omega}(v_0)\bigr).
\]
As the reference values are equal, no function
$\widehat m:[0,\infty)\to\mathbb R$ can satisfy
\[
 m\bigl(\widetilde{\mathcal H}_{\eta,\omega}(u)\bigr)
 =\widehat m\bigl(\mathcal H(u)\bigr)
 \quad\text{for every }u\in\mathcal D.
\]
The obstruction already occurs for strictly positive functions in
$\mathcal D$, since both $v_0$ and $v_\varepsilon$ are positive in $\Omega$.

If $|U|=|\Omega|$, then $J=\mathcal H$, and the perturbation can be
absorbed into $m$ by setting
\[
 \widehat m(s)=m\bigl(s[1+\eta\sin(\omega s)]\bigr),
 \qquad s\geq0.
\]
\end{example}

\section*{Acknowledgements}
The first author is a member of the Gruppo Nazionale per l'Analisi Matematica, la Probabilit\`a e le loro Applicazioni (GNAMPA) of the Istituto Nazionale di Alta Matematica (INdAM). 
This work was partially supported by the INdAM GNAMPA Project
with CUP E5324001950001.

\section*{Conflict of interest}
The authors declare that they have no conflict of interest.

\end{document}